\documentclass[USenglish]{article}

\usepackage[T1]{fontenc}
\usepackage[utf8]{inputenc} 
\usepackage{microtype}
\usepackage{amssymb}
\usepackage{amsmath}
\usepackage{graphicx}
\usepackage{authblk}

\newcommand{\mb}[1]{\mathbf{#1}}
\newcommand{\bs}[1]{\boldsymbol{#1}}

\newcommand{\cP}{\mathcal{P}}

\newcommand{\cF}{\mathcal{F}}
\newcommand{\bR}{\mathbb{R}}
\newcommand{\bC}{\mathbb{C}}
\newcommand{\cS}{\mathcal{S}}
\newcommand{\tp}{^{\mbox{\tiny\bf T}}}

\newcommand{\proof}{\textbf{Proof}\hspace{0.25cm}}
\newcommand{\qed}{\hfill $\square$}
\newcommand{\keywords}{\textbf{Keywords}\hspace{0.25cm}}

\newtheorem{problem}{Problem}[section]
\newtheorem{corollary}{Corollary}[section]
\newtheorem{theorem}{Theorem}[section]

\begin{document}

\title{The modulus of the Fourier transform on a sphere determines $3$-dimensional convex polytopes}

\author{Konrad Engel\thanks{konrad.engel@uni-rostock.de} }
\author{Bastian Laasch\thanks{bastian.laasch@uni-rostock.de}}

\affil{University of Rostock, Institute for Mathematics, 18057 Rostock, Germany}

\maketitle

% \runningauthor{Konrad Engel and Bastian Laasch}
% \runningtitle{The modulus of the FT determines $3$-dim polytopes}

\abstract{Let $\cP$ and $\cP'$ be $3$-dimensional convex polytopes in $\bR^3$ and $S \subseteq \bR^3$ be a non-empty intersection of an open set with a sphere.
As a consequence of a somewhat more general result it is proved that $\cP$ and $\cP'$ coincide up to translation and/or reflection in a point if $|\int_{\cP} e^{-i\mb{s}\cdot\mb{x}} \,\mb{dx}| = |\int_{\cP'} e^{-i\mb{s}\cdot\mb{x}} \,\mb{dx}|$ for all $\mb{s} \in S$. This can be applied to the field of crystallography regarding the question whether a nanoparticle modelled as a convex polytope is uniquely determined by the intensities of its X-ray diffraction pattern on the Ewald sphere.\\}
  
  
% \classification[MSC]{42B10; 52B10; 52B11; 81U40}

\section{Introduction}\label{sec::1}

In this paper, we prove the statement of the title: The modulus of the Fourier transform on a sphere determines $3$-dimensional convex polytopes.
In small-angle X-ray scattering and partly also in wide-angle X-ray scattering of nanoparticles the modulus of the Fourier transform of the reflected beam wave vectors can be ``approximately'' measured on the \emph{Ewald sphere}, see e.g. \cite{Raines2010}, \cite{Barke2015}, \cite{Rossbach2019}. Hence, uniqueness questions are of special interest and by our result, at least theoretically, the polytopes can be reconstructed from the measurements.

The statement can be formulated in a more general form and therefore we need several definitions and notations first.

The \emph{Fourier transform} $\cF_f$ of an $L_1$-integrable function $f:\bR^n\rightarrow\bR$ is defined by
\begin{align*}
    \cF_f(\mb{s})=\int_{\bR^n} f(\mb{x})e^{-i\mb{s}\cdot\mb{x}} \,\mb{dx}\,,\quad\mb{s}\in\bR^n\,,
\end{align*}
where the product $\cdot$ is the standard scalar product, see e.g. \cite{Plonka2018}. If $f=\chi_{\cP}$ is the characteristic function of an $n$-dimensional compact set $\cP\subseteq\bR^n$, then
\begin{align*}
    \cF_f(\mb{s})=\int_{\cP} e^{-i\mb{s}\cdot\mb{x}} \,\mb{dx}
\end{align*}
is the Fourier transform of $\cP$. We denote it briefly by $\cF_{\cP}(\mb{s})$.

An important question is the following: Assume that $\cP$ and $\cP'$ are $n$-dimensional compact sets in $\bR^n$ and that 
$ |\cF_{\cP}(\mb{s})|=|\cF_{\cP'}(\mb{s})|$
for all $\mb{s}$ belonging to a sparse subset of $\bR^n$. Is it true that $\cP$ and $\cP'$ have the same ``structure''?

We say that $\cP$ and $\cP'$ are \emph{strongly congruent}, denoted by $\cP \cong_s \cP'$, if there is some vector $\mb{v}$ and an $\epsilon \in \{-1,1\}$ such that
\begin{align}
\label{strongly congruent}
    \cP'=\epsilon \cP + \mb{v}\,,
\end{align}
i.e., $\cP$ and $\cP'$ coincide up to \emph{translation} and/or \emph{reflection in a point}.
One easily obtains from (\ref{strongly congruent}) that for all $\mb{s} \in \bR^n$
\begin{align*}
    \cF_{\cP'}(\mb{s}) = 
    \begin{cases}
        e^{i \mb{s} \cdot \bs{v}} \cF_{\cP}(\mb{s})&\text{ if } \epsilon=1\,,\\
        e^{i \mb{s} \cdot \bs{v}} \overline{\cF_{\cP}(\mb{s})}&\text{ if } \epsilon=-1\,,
    \end{cases}
\end{align*}
which implies
\begin{align*}
    |\cF_{\cP'}(\mb{s})| = |\cF_{\cP}(\mb{s})|\,.
\end{align*}

Hence, two strongly congruent sets cannot be distinguished by the modulus of their Fourier transform and the above question has to be formulated more accurately, e.g. in form of a problem as follows:

\begin{problem}
\label{FT-problem}
Determine sufficient conditions for $\cP, \cP', S \subseteq \bR^n$ such that 
\begin{align}
\label{FT-problem equation}
    \forall \mb{s} \in S \quad |\cF_{\cP}(\mb{s})| = |\cF_{\cP'}(\mb{s})| \Longrightarrow \cP\cong_s\cP'\,.
\end{align}
\end{problem}

We will see that in the case $S=\bR^n$ this problem is equivalent to the so-called \emph{covariogram problem}, introduced by Matheron in \cite{Matheron1986}. 

The \emph{covariogram} of an $n$-dimensional compact set $\cP\subseteq\bR^n$ is a function $g_{\cP}: \bR^n \rightarrow \bR$ defined by
\begin{align*}
    g_{\cP}(\mb{x})=\lambda_n (\cP\cap(\cP+\mb{x}))\,,\quad\mb{x}\in\bR^n\,,
\end{align*}
where $\lambda_n$ denotes the $n$-dimensional Lebesgue measure, see \cite{Matheron1975}. Thus $g_{\cP}$ associates with each $\mb{x}$ the volume of the intersection of $\cP$ with the \emph{Minkowski sum} $\cP+\mb{x}$.
It is again easy to see that $\cP \cong_s \cP'$ implies $g_{\cP}(\mb{x})=g_{\cP'}(\mb{x})$ for all $\mb{x} \in \bR^n$, i.e., the covariogram is invariant with respect to translation and reflection in a point.

\begin{problem}[covariogram problem]
\label{covariogram problem}
Determine sufficient conditions for $\cP,$ $\cP'$ $\subseteq \bR^n$ such that 
\begin{align}
    \label{covariogram problem equation}
    \forall \mb{x} \in \bR^n \quad g_{\cP}(\mb{x}) = g_{\cP'}(\mb{x}) \Longrightarrow \cP\cong_s\cP'\,.
\end{align}
\end{problem}

The equivalence of both problems (see e.g. \cite{Bianchi2009}) follows from the following theorem which we reprove in Section \ref{sec::2} in order to make the paper self-contained.

\begin{theorem}
\label{equivalent problems}
We have
\begin{align}
\label{equivalent problems equation}
    \forall \mb{s} \in \bR^n \quad |\cF_{\cP}(\mb{s})| = |\cF_{\cP'}(\mb{s})| \Longleftrightarrow \forall \mb{x} \in \bR^n \quad g_{\cP}(\mb{x}) = g_{\cP'}(\mb{x})\,.   
\end{align}
\end{theorem}

Fortunately, there exist deep results in the literature concerning Problem \ref{covariogram problem}. The following conditions are sufficient for (\ref{covariogram problem equation}) (see also Theorem 1 of the survey paper \cite{Horvath2020}):
\begin{itemize}
    \item $\cP, \cP' \subseteq \bR^2$ are $2$-dimensional convex polygons \cite[Theorem 3.1]{Nagel1993}.
    \item $\cP, \cP' \subseteq \bR^2$ are $2$-dimensional convex and compact sets \cite[Theorem 1.1]{Averkov2009}.
    \item $\cP, \cP' \subseteq \bR^3$ are $3$-dimensional convex polytopes  \cite[Theorem 1.1]{Bianchi2009}.
    \item $\cP, \cP' \subseteq \bR^n$ are $n$-dimensional ($n\geq 3$) convex simplicial polytopes (i.e., each facet is a simplex and each polytope is in general relative position to its reflection, see \cite[Corollary to Theorem 2]{Goodey1997} and \cite{Schneider1998}). Note that, for example, in the $3$-dimensional case this is fulfilled for a regular tetrahedron but not for a regular icosahedron.
\end{itemize}

In other words, these bodies are determined by their covariogram, respectively by the modulus of their Fourier transform, over $\bR^n$ up to translation and reflection in a point. Note that for arbitrary convex sets in $\bR^3$ the question is still open and for dimension $n\geq4$ counterexamples can be constructed, e.g. for special convex polytopes (see \cite[Theorem 1.2]{Bianchi2005} and \cite[Remark 9.3]{Bianchi2009}).

The more general problem of the reconstruction of an arbitrary function $f$ from the modulus of its Fourier transform over $\bR^n$ is called \emph{phase retrieval problem} and was intensively studied in the past (see e.g. \cite{Hurt1989, Klibanov1995, Rosenblatt1984} and the references given in \cite{Bianchi2002}). Usually, phase retrieval is under-determined without any additional constraints. Hence, a priori assumptions that the function has a particular form, e.g. $f$ is the characteristic function of an $n$-dimensional polytope, are needed.

The problem will become relevant for physical applications if we assume the identity $|\cF_{\cP}(\mb{s})| = |\cF_{\cP'}(\mb{s})|$ not over $\bR^n$, like in (\ref{equivalent problems equation}), but only over a subset $S \subseteq \bR^n$, like in (\ref{FT-problem equation}). We can use the above listed results to solve Problem \ref{FT-problem} if we have a solution of the following one:

\begin{problem}
\label{sparse FT-problem}
Determine sufficient conditions for $\cP, \cP', S \subseteq \bR^n$ such that 
\begin{align}
\label{sparse FT-problem equation}
    \forall \mb{s} \in S \quad |\cF_{\cP}(\mb{s})| = |\cF_{\cP'}(\mb{s})| \Longrightarrow \forall \mb{s} \in \bR^n \quad |\cF_{\cP}(\mb{s})| = |\cF_{\cP'}(\mb{s})| \,.
\end{align}
\end{problem}

Section \ref{sec::3} is devoted to this problem. In Section \ref{sec::4} we discuss the physical application for X-ray scattering and, in particular, derive the following corollary as an answer to the title of this paper:

\begin{corollary}
\label{sphere}
If $\cP$ and $\cP'$ are 3-dimensional convex polytopes in $\bR^3$ and $S \subseteq \bR^3$ is a non-empty intersection of an open set with a sphere, then 
%(\ref{sparse FT-problem equation}) is true and consequently
\[
\forall \mb{s} \in S \quad |\cF_{\cP}(\mb{s})| = |\cF_{\cP'}(\mb{s})| \Longrightarrow \cP\cong_s\cP'\,.
\]
%thus also (\ref{FT-problem equation}) are true.
\end{corollary}
Finally, Section \ref{sec::5} contains a generalization to a further class of $n$-dimensional convex polytopes in addition to the already mentioned simplicial polytopes.

\section{Proof of Theorem \ref{equivalent problems}}\label{sec::2}

Recall that the convolution of two functions $f,g:\bR^n\rightarrow\bR$ with regard to the measure $\lambda_n$ is defined by 
\begin{align*}
    (f\ast g)(\mb{x})=\int_{\bR^n} f(\mb{y})g(\mb{x}-\mb{y}) \,\mb{d}\lambda_n(\mb{y})\,,\quad\mb{x}\in\bR^n\,.
\end{align*}
Hence, we have
\begin{align*}
    (\chi_{\cP}\ast\chi_{-\cP})(\mb{x})
    &=\int_{\bR^n} \chi_{\cP}(\mb{y})\chi_{-\cP}(\mb{x}-\mb{y}) \,\mb{d}\lambda_n(\mb{y})\\
    &=\int_{\bR^n} \chi_{\cP}(\mb{y})\chi_{\cP}(\mb{y}-\mb{x}) \,\mb{d}\lambda_n(\mb{y})\\
		&=\int_{\bR^n} \chi_{\cP}(\mb{y})\chi_{\cP+\mb{x}}(\mb{y}) \,\mb{d}\lambda_n(\mb{y})\\
    &=\int_{\cP\cap(\cP+\mb{x})} 1 \,\mb{d}\lambda_n(\mb{y})
    =\lambda_n (\cP\cap(\cP+\mb{x}))\,,
\end{align*}
and consequently
\begin{align}
\label{covariogram and convolution}
    g_{\cP}(\mb{x})=(\chi_{\cP}\ast\chi_{-\cP})(\mb{x})\quad\forall\mb{x}\in\bR^n\,.
\end{align}

Now well-known properties of the Fourier transform imply that the following statements are equivalent:
\begin{align*}
    g_{\cP}(\mb{x}) = g_{\cP'}(\mb{x}) \quad \forall \mb{x} \in \bR^n\,
    &\Longleftrightarrow \cF_{g_{\cP}}(\mb{s})=\cF_{g_{\cP'}}(\mb{s}) \quad &\forall \mb{s} \in \bR^n\,~\\
    &\Longleftrightarrow \cF_{\cP}(\mb{s}) \cF_{-\cP}(\mb{s})=\cF_{\cP'}(\mb{s}) \cF_{-\cP'}(\mb{s}) \quad &\forall \mb{s} \in \bR^n\,~\\
     &\Longleftrightarrow \cF_{\cP}(\mb{s}) \overline{\cF_{\cP}(\mb{s})}=\cF_{\cP'}(\mb{s}) \overline{\cF_{\cP'}(\mb{s})} \quad &\forall \mb{s} \in \bR^n\,~\\
     &\Longleftrightarrow |\cF_{\cP}(\mb{s})|=|\cF_{\cP'}(\mb{s})|\quad &\forall \mb{s} \in \bR^n\,.
\end{align*}
The first equivalence holds because of the injectivity of the Fourier transform (see \cite{Plonka2018}). The second one follows from (\ref{covariogram and convolution}) and the identity between the Fourier transform of two convoluted functions with the product of their Fourier transforms. For the last equivalence we used $|z|^2=z\overline{z}$ for $z\in\bC$.
\qed

\section{An answer to Problem \ref{sparse FT-problem}}\label{sec::3}

In this paper, we not only allow the restriction of the Fourier transform to a sphere, but more generally to rationally parameterisable hypersurfaces satisfying two general conditions. Therefore, we need some definitions from \cite{Engel2020}.

A \emph{rationally parameterisable hypersurface} (briefly \emph{rp-hypersurface}) is a set $\cS$ of points in $\bR^n$ of the form
\begin{align}\label{rps}
    \cS=\{\bs{\sigma}(\mb{t}): \mb{t} \in D\}\,,
\end{align}
where 
\begin{align*}
    \bs{\sigma}(\mb{t})=
    \begin{pmatrix}
        \sigma_1(t_1,\dots,t_{n-1})\\
        \vdots\\
        \sigma_n(t_1,\dots,t_{n-1})
    \end{pmatrix}\,,
\end{align*}
the functions $\sigma_j$ are rational functions, $j=1,\dots,n$, and $D \subseteq \bR^{n-1}$ is the domain of $\cS$. 

Using spherical coordinates and the standard substitution $t=\tan(\frac{\alpha}{2})$, which implies $\cos(\alpha)=\frac{1-t^2}{1+t^2}$ and $\sin(\alpha)=\frac{2t}{1+t^2}$, one obtains that the $3$-dimensional unit sphere 
\begin{align*}
    \left\{\left(\frac{2t_1}{1+t_1^2}\frac{1-t_2^2}{1+t_2^2},\frac{2t_1}{1+t_1^2}\frac{2t_2}{1+t_2^2},\frac{1-t_1^2}{1+t_1^2}\right): t_1\in \bR_+, t_2\in\bR\right\}
\end{align*}
with the missing segment $(-\sqrt{1-\lambda^2},0,\lambda), \lambda\in[-1,1)$, is an rp-hypersurface. Since affine transformations do not violate the rationality this is also true for any sphere in $\bR^3$.

With the function $\bs{\sigma}: \bR^{n-1} \rightarrow \bR^n$ we associate its \emph{normalized function}, i.e., the function $\hat{\bs{\sigma}}: \bR^{n-1} \rightarrow \bR^{n-1}$ defined by
\begin{align*}
    \hat{\bs{\sigma}}(\mb{t})=
    \begin{pmatrix}
        \frac{\sigma_2(t_1,\dots,t_{n-1})}{\sigma_1(t_1,\dots,t_{n-1})}\\
        \vdots\\
        \frac{\sigma_n(t_1,\dots,t_{n-1})}{\sigma_1(t_1,\dots,t_{n-1})}
    \end{pmatrix}\,.
\end{align*}
Note that $\hat{\bs{\sigma}}(\mb{t})$ is only defined if $\mb{t} \in D \setminus\sigma_1^{-1}(0)$, where $\sigma_1^{-1}(0)$ is the zero set of $\sigma_1$.
For a set $O \subseteq D\setminus\sigma_1^{-1}(0)$ let 
\begin{align*}
    \hat{\bs{\sigma}}(O)=\{\hat{\bs{\sigma}}(\mb{t}): \mb{t} \in O\}\,.
\end{align*}

In the following we need that, for an rp-hypersurface $\cS$ given by (\ref{rps}), an open subset $O$ of $D\setminus\sigma_1^{-1}(0)$ in $\bR^{n-1}$ satisfies two conditions:

\textbf{Hyperplane condition}: $\bs{\sigma}(O)$ is not contained in a hyperplane.

\textbf{Inner point condition}: There is a $\mb{t} \in O$ such that $\hat{\bs{\sigma}}(\mb{t})$ is an inner point of $\hat{\bs{\sigma}}(O)$ in $\bR^{n-1}$.

The proof of the following Theorem \ref{rp-hyperface theorem} is based on a result on E-functions of \cite[Theorem 2.3 together with Theorem 2.2. and Lemma 2.2]{Engel2020}. 
A function $F:\bR^n\rightarrow\bC$ is called an \emph{E-function of degree $d$} if it has the form
\begin{align*}
    F(\mb{s})=\sum_{k=1}^m P_k(\mb{s})e^{-i\mb{v}_k\cdot\mb{s}}
\end{align*}
with distinct points $\mb{v}_k$, $k=1,\dots,m$, and homogeneous rational functions $P_k(\mb{s})$ of degree $d$.

\begin{theorem}[see \cite{Engel2020}]  
\label{E-functions theorem}
Let $\cS=\{\bs{\sigma}(\mb{t}): \mb{t} \in D\}$ be an rp-hypersurface and let  $O \subseteq D\setminus\sigma_1^{-1}(0)$ be an open subset of $\bR^{n-1}$ that  satisfies the hyperplane and the inner point condition.
Let $F(\mb{s})$ be an E-function of degree $d$. If 
\begin{align*}
    F(\mb{s}) = 0\, \quad \forall \mb{s} \in \bs{\sigma}(O)\,,
\end{align*}
then all coefficients of the exponential functions and hence also $F$ are the zero function (up to the cases where some denominator equals 0).
\end{theorem}

Now we are able to present an answer to Problem \ref{sparse FT-problem}.

\begin{theorem}
\label{rp-hyperface theorem}
Let $\cS=\{\bs{\sigma}(\mb{t}): \mb{t} \in D\}$ be an rp-hypersurface in $\bR^n$ and let $O\subseteq D\setminus\sigma_1^{-1}(0)$  be an open subset of $\bR^{n-1}$ that satisfies the hyperplane and the inner point condition.
If $\cP$ and $\cP'$ are $n$-dimensional convex polytopes in $\bR^n$ and $S=\bs{\sigma}(O)$, then (\ref{sparse FT-problem equation}) %and thus also (\ref{FT-problem equation}) are 
is true.
\end{theorem}

\proof
For an $n$-dimensional convex polytope $\cP$ let $V_{\cP}\subseteq\bR^n $ be its vertex set and  let $E_{\cP}=\{\mb{v}_1-\mb{v}_2:\mb{v}_1,\mb{v}_2\in V_{\cP}\}$. For the Fourier transform of a convex polytope $\cP$ it is known (see \cite{Barvinok1994, Brion1988, Davis1964, Lawrence1991, Pukhlikov1992}) that
\begin{align}
\label{ft of convex polytope}
    \cF_{\cP}(\mb{s})=\sum_{\mb{v}\in V_{\cP}} Q_{\cP, \mb{v}}(\mb{s}) e^{-i\mb{v}\cdot\mb{s}}\,\quad\forall\mb{s}\in\bR^n\setminus Z_{\cP}\,,
\end{align}
where each coefficient is a rational function of the form
\begin{align*}
    Q_{\cP, \mb{v}}(\mb{s})=\sum_{I\in\mathcal{I}_{\cP}} \frac{\lambda_{\cP, \mb{v}, I}}{\prod_{\mb{e}\in I}\mb{e}\cdot\mb{s}}\,,
\end{align*}
which is not the zero function, see \cite[Lemma 2]{Gravin2012596}.
The numerators $\lambda_{\cP, \mb{v}, I}$ are real numbers and $\mathcal{I}_{\cP}$ is a family of $n$-element linearly independent subsets of $E_{\cP}$. The set $Z_{\cP}$ contains only vectors $\mb{s}$ for which the scalar product $\mb{e}\cdot\mb{s}$ vanishes. Therefore, (\ref{ft of convex polytope}) is an E-function of degree $-n$ and since $|z|^2=z\overline{z}$ for $z\in\bC$ we have for the squared modulus of the Fourier transform
\begin{align}
\label{squared ft}
    |\cF_{\cP}(\mb{s})|^2=\left(\sum_{\mb{v}\in V_{\cP}} Q_{\cP, \mb{v}}(\mb{s})^2\right) e^{-i\mb{0}\cdot\mb{s}} + \sum_{\substack{\mb{v}_i,\mb{v}_j\in V_{\cP}\\\mb{v}_i\neq\mb{v}_j}} Q_{\cP, \mb{v}_i}(\mb{s})Q_{\cP, \mb{v}_j}(\mb{s})e^{-i(\mb{v}_i-\mb{v}_j)\cdot\mb{s}}\,.
\end{align}
It is possible that the exponents of the exponential functions in the second sum of (\ref{squared ft}) are not distinct. In such cases we can merge the corresponding terms and see that (\ref{squared ft}) is an E-function of degree $-2n$. Note that a linear combination of E-functions of degree $d$ is again an E-function of degree $d$. Therefore, $|\cF_{\cP}(\mb{s})|^2-|\cF_{\cP'}(\mb{s})|^2$ is also an E-function of degree $-2n$.

The assumption in (\ref{sparse FT-problem equation}), respectively (\ref{FT-problem equation}), implies
\begin{align}
\label{squared ft difference}
    |\cF_{\cP}(\mb{s})|^2-|\cF_{\cP'}(\mb{s})|^2=0\,\quad \forall \mb{s} \in \bs{\sigma}(O)\
\end{align}
and by Theorem \ref{E-functions theorem} and the continuity of the Fourier transform in $\mb{s}$
\begin{align*}
    |\cF_{\cP}(\mb{s})|-|\cF_{\cP'}(\mb{s})|=0\,\quad \forall \mb{s} \in \bR^n\,.
\end{align*}
\qed

\section{Proof of Corollary \ref{sphere} and a physical application}\label{sec::4}

\proof 
In \cite{Engel2020} the following theorem is proved for quadratic hypersurfaces, i.e., 
zero sets of an equation of the form
\begin{align*}
    \frac{1}{2}\mb{s}\tp A \mb{s} + \mb{b}\tp \mb{s} + c = 0\,,
\end{align*}
where $A$ is a symmetric matrix different from the zero matrix.

\begin{theorem}
\label{quadric hypersurface theorem}
If $\cS$ is a quadratic hypersurface that does not contain a line but at least two points, then, up to an exceptional set of  hypersurface measure zero, it is an rp-hypersurface with some parameter domain $D$ and every open subset $O$ of $D\setminus\sigma_1^{-1}(0)$ in $\bR^{n-1}$ satisfies the hyperplane and the inner point condition.
\end{theorem}

Clearly, this theorem can be applied to a sphere since a sphere is a quadratic hypersurface that does not contain a line but at least two points.

We may assume that the sphere is the unit sphere and that the non-empty intersection $S$ of the open set with the sphere does not contain points from the missing segment $(-\sqrt{1-\lambda^2},0,\lambda), \lambda\in[-1,1)$, in the parametrization of the unit sphere and no points with vanishing first coordinate. Obviously, the set 
\begin{align*}
    O=\left\{(t_1,t_2): \left(\frac{2t_1}{1+t_1^2}\frac{1-t_2^2}{1+t_2^2},\frac{2t_1}{1+t_1^2}\frac{2t_2}{1+t_2^2},\frac{1-t_1^2}{1+t_1^2}\right) \in S\right\}
\end{align*}
has the properties required in Theorem \ref{quadric hypersurface theorem} and thus Corollary \ref{sphere} follows from Theorem \ref{rp-hyperface theorem}, Theorem \ref{equivalent problems} and the fact that 3-dimensional convex polytopes are positive examples for Problem \ref{covariogram problem}. 
\qed\newline

Now we apply our results to the field of crystallography. An interesting questions is whether the $3$-dimensional structure of nanoparticles is determined by their single-shot X-ray diffraction pattern (see \cite{Raines2010}, \cite{Barke2015} and, for a current list of references, \cite{Rossbach2019}). 

\begin{figure}[htpb]
    \centering
    \includegraphics[width=10cm]{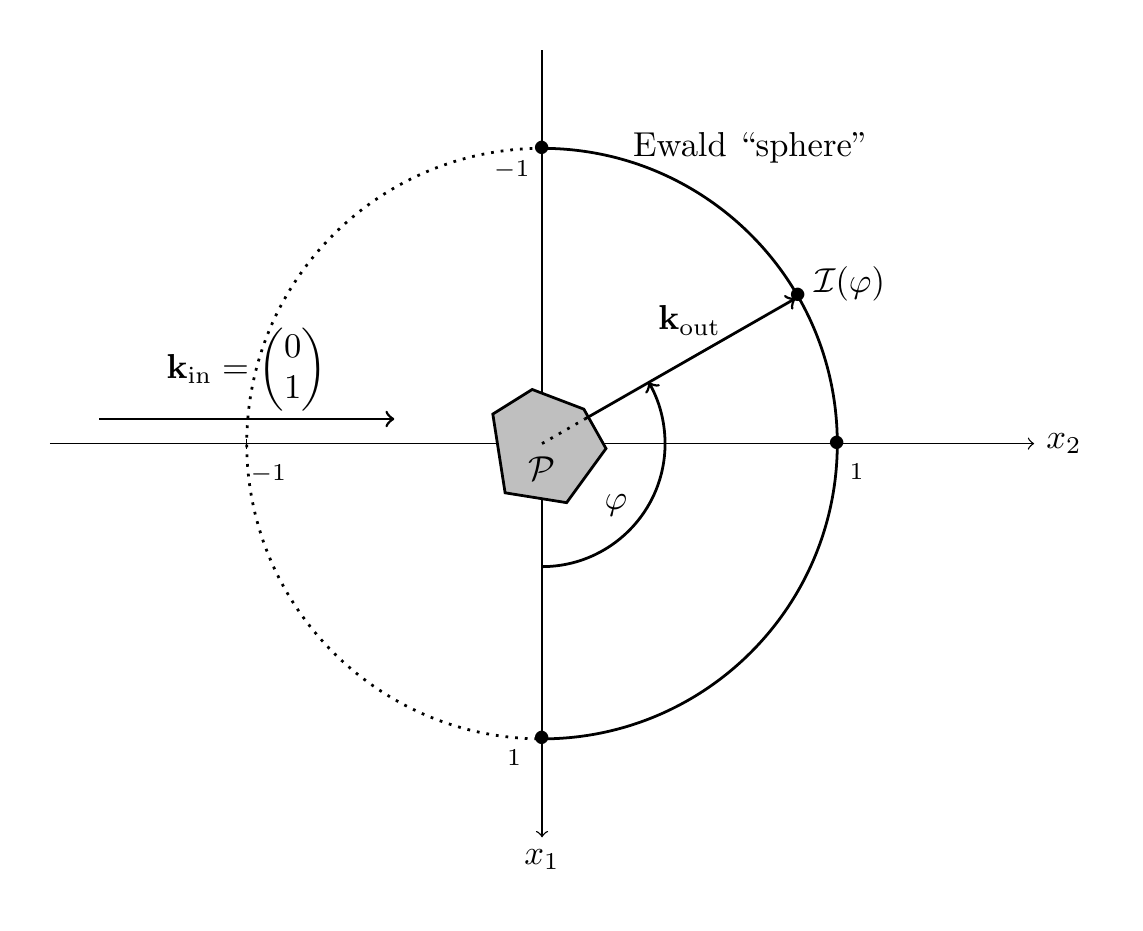}
    \caption{The $2$-dimensional single-shot scattering model for a convex polygon, respectively nanoparticle, $\cP$ with normalized wavelengths for the incoming and reflected X-ray beam wave vectors $\textbf{k}_\text{in}$ and $\textbf{k}_\text{out}$. The solid arc marks the area of the ($2$-dimensional) Ewald ``sphere'' where the scattering intensities $\mathcal{I}$ can be measured during experiments.}
    \label{scattering model}
\end{figure}

To this end we consider the following $2$-dimensional setup (see Figure \ref{scattering model}) -- the $3$-dimensional case follows analogously. 
The incoming normalized X-ray beam wave vector $\mb{k}_\text{in}=(0,1)^T$ with wavelength $1$ illuminates a nanoparticle $\cP$ modelled as a convex polygon. The occurring diffraction event causes scattered wave vectors $\mb{k}_\text{out}$ in different directions. The set of all such vectors constitutes the so-called diffraction pattern. We assume that during the diffraction process there is no energy gained or lost. Hence,
\begin{align*}
    \| \mb{k}_\text{out}\|_2=\| \mb{k}_\text{in}\|_2=1\,,
\end{align*}
i.e., the wavelengths of the diffracted beams are also $1$. Therefore, in the $2$-dimensional case the vectors $\mb{k}_\text{out}$ form a unit circle, respectively in the $3$-dimensional case a unit sphere, the Ewald sphere. Note that we assume that the diffracted wave vectors $\mb{k}_\text{out}$ are originated in the center of the Ewald sphere since the particle $\cP$ is infinitesimally small. The intensity $\mathcal{I}$ of a scattered wave is proportional to
\begin{align}
\label{scattering intensity}
    \left|\int_{\cP} e^{-i(\mb{k}_\text{out}-\mb{k}_\text{in})\cdot\mb{x}} \,\mb{dx}\right|\,,
\end{align}
the difference $\mb{k}_\text{out}-\mb{k}_\text{in}$ of the reflected and the incoming wave vector is called scattering vector, i.e, the intensity is proportional to the modulus of the Fourier transform of the scattering vector. 

In experiments we can measure these intensities on a semicircle of the Ewald ``sphere'' with exception of the points $(1,0)^T$, $(0,1)^T$ and $(-1,0)^T$ (see the solid arc in Figure \ref{scattering model}). Representing $\mb{k}_\text{out}$ in polar coordinates yields
\begin{align*}
    \mb{k}_\text{out}=\begin{pmatrix} \cos(\varphi)\\ \sin(\varphi)\\ \end{pmatrix}\,, \quad \varphi\in (0,\pi)\setminus\left\{\frac{\pi}{2}\right\}\,,
\end{align*}
and together with (\ref{scattering intensity}) we get
\begin{equation}
\label{scattering intensity 2}
    \begin{aligned}
        \mathcal{I}(\varphi)&\propto
        \left|\int_{\cP} e^{-i(\cos(\varphi)x_1+(\sin(\varphi)-1)x_2)} \,\mb{dx}\right|\\
        &=\left|\cF_{\cP}\left(\begin{pmatrix} \cos(\varphi)\\ \sin(\varphi)-1\\ \end{pmatrix}\right)\right|\,.
    \end{aligned}
\end{equation}
Hence, the set $S$ mentioned in the assumptions of Problem \ref{FT-problem} and Problem \ref{sparse FT-problem} is given by
\begin{align*}
    S=\left\{\begin{pmatrix} \cos(\varphi)\\ \sin(\varphi)-1\\ \end{pmatrix}:\varphi\in (0,\pi)\setminus\left\{\frac{\pi}{2}\right\}\right\}\,,
\end{align*}
i.e., a semicircle shifted in the second coordinate.

Analogously to (\ref{scattering intensity 2}) we can compute the intensities for the $3$-dimensional case on a set $S$ corresponding to a (translated) hemisphere of the Ewald sphere with some excluded points in spherical coordinates
\begin{equation*}
    \begin{aligned}
        \mathcal{I}(\varphi,\theta)&\propto
        \left|\int_{\cP} e^{-i(\sin(\theta)\cos(\varphi)x_1+\sin(\theta)\sin(\varphi)x_2+(\cos(\theta)-1)x_3)} \,\mb{dx}\right|\\
        &=\left|\cF_{\cP}\left(\begin{pmatrix} \sin(\theta)\cos(\varphi)\\ \sin(\theta)\sin(\varphi)\\\cos(\theta)-1\\ \end{pmatrix}\right)\right|\,.
    \end{aligned}
\end{equation*}
By Corollary \ref{sphere} we know that two $3$-dimensional convex polytopes with the same diffraction pattern on the hemisphere are strongly congruent. So the measured intensities on the Ewald sphere of scattered X-ray beams determine the underlying object uniquely up to translation and reflection in a point. Note that the experimental setup delivers ``only'' finitely many approximated intensities and we assumed infinitely many exact measurements.

\section{A generalization to a further class of $n$-dimensional convex polytopes}\label{sec::5}

If $\cP$ and $\cP'$ are convex polytopes of arbitrary dimension with vertex sets $V_{\cP}$ and $V_{\cP'}$, then clearly
\begin{equation}
\label{congruent vertex sets}
    V_{\cP} \cong_s V_{\cP'} \Longrightarrow \cP\cong_s\cP'\,.
\end{equation}
Thus, we are lead to a discrete variant.

The points of a finite set $V$ in $\bR^n$ are in \emph{general position} if each line through the origin contains, in addition to the origin, at most one pair $\pm\mb{w}$ of the multiset $V-V$. The following theorem is proved in \cite[Theorem 1]{Senechal2008}:
\begin{theorem}
\label{general position sets}
    Let $V$ and $V'$ be finite subsets of $\bR^n$, each with points in general position. Then
    \[
        V-V=V'-V'\Longrightarrow V\cong_s V'\,.
    \]
\end{theorem}
This theorem enables us to present another answer to Problem \ref{FT-problem} for the $n$-dimensional case.

\begin{theorem}
\label{general position polytopes}
    Let $\cS=\{\bs{\sigma}(\mb{t}): \mb{t} \in D\}$ be an rp-hypersurface in $\bR^n$ and let $O\subseteq D\setminus\sigma_1^{-1}(0)$  be an open subset of $\bR^{n-1}$ that satisfies the hyperplane and the inner point condition.
    Further let $\cP$ and $\cP'$ be $n$-dimensional convex polytopes in $\bR^n$ such that the vertices of $\cP$ as well of $\cP'$ are in general position.
    Finally let $S=\bs{\sigma}(O)$. Then
    \[
    \forall \mb{s} \in S \quad |\cF_{\cP}(\mb{s})| = |\cF_{\cP'}(\mb{s})| \Longrightarrow \cP\cong_s\cP'\,.
    \]
\end{theorem}
\proof 
We use again (\ref{squared ft}) and (\ref{squared ft difference}). Since the coefficients of the exponential functions are not the zero function we obtain that $V_{\cP}-V_{\cP}=V_{\cP'}-V_{\cP'}$. Now the assertion follows immediately from Theorem \ref{general position sets} and (\ref{congruent vertex sets}). \qed\newline

Note that this theorem can be applied in particular to simplices because their vertex sets are in general position.

\section{Conclusion}\label{sec::6}

In this paper we used a series of results of discrete geometry, namely solutions for the covariogram problem, to answer the question whether the intensities of an X-ray diffraction pattern on the Ewald sphere determine the underlying object uniquely. For physical application the $3$-dimensional case is of interest and if the illuminated nanoparticle is modelled as a convex polytope the reconstruction is unique up to translation and reflection in a point. However, our theoretically result assumes infinitely many exact measurements what is not the case in experimental setups. Therefore, the question remains open whether the result is true under the condition of finitely many given intensities.

\section*{Acknowledgement}
We would like to thank Thomas Fennel and Stefan Scheel for presenting the motivation of this study.

This work was partly supported by the European Social Fund (ESF) and the Ministry of Education, Science and Culture of Mecklenburg-Western Pomerania (Germany) within the project NEISS – Neural Extraction of Information, Structure and Symmetry in Images under grant no ESF/14-BM-A55-0006/19.

\bibliography{references}
\bibliographystyle{plain}

\end{document}